\documentclass{amsart}

\usepackage{amssymb,latexsym}

\begin{document}

\title{\bf A counterexample to the "composition conjecture".}
\author {F. Pakovich}
\address{Department of Mathematics, Weizmann Institute of Science,
Rehovot 76100, Israel}
\email{pakovich@wisdom.weizmann.ac.il}

\subjclass{Primary 34C99; Secondary 30D05}

\date{\today}

\keywords{Poincare center-focus problem; polynomial Abel equation;
polynomials; compositions; functional equations}

\begin{abstract}
In this note we construct a class of counterexamples to the "composition
conjecture" concerning an infinitesimal version of the center problem for
the polynomial Abel equation in the complex domain.
\end{abstract}

\maketitle

\def\d{{\rm d}}
\def\D{{\rm D}}
\def\I{{\rm I}}

\def\C{{\mathbb C}}
\def\N{{\mathbb N}}
\def\P{{\mathbb P}}   
\def\Z{{\mathbb Z}}
\def\d{{\rm d\,}}
\def\deg{{\rm deg\,}}
\def\Det{{\rm Det}}\def\dim{{\rm dim\,}}
\def\Ker{{\rm Ker\,}}
\def\Gal{{\rm Gal\,}}
\def\St{{\rm St\,}}
\def\Sym{{\rm Sym\,}}
\def\Mon{{\rm Mon\,}}

\vskip 0.2cm
In this note we treat the following "polynomial moment problem"
proposed in \cite{c1}, \cite{c2} as an infinitesimal version of
the center problem for the polynomial Abel equation in the complex domain:
{\it for complex polynomials $P(z),Q(z)$ to find
conditions under which all moments $$m_{i}(P,Q,a,b)=\int^b_a
P^i(z)\d Q(z)=\int^b_a P^i(z)Q^{\prime}(z)\d z, \ \ i\geq 1, \ \ a,b\in 
\C,$$ vanish under assumption that $P(a)=P(b),$ $Q(a)=Q(b),$ $a\neq b.$ }

The "composition conjecture" suggested in \cite{c1}
states that the following condition is necessary and
sufficient:
{\it there exist polynomials $\tilde P(z),
\tilde Q(z), W(z)$ such that $$P(z)=\tilde P(W(z)), \ \ \ \
Q(z)=\tilde Q(W(z)),\ \ \ \ \deg W(z)>1, \eqno(*) $$ and $W(a)=W(b).$}
Note that as $m_{i}(P,Q,a,b)=m_{i}(\tilde P,\tilde Q,W(a),W(b))$ and
$W(a)=W(b)$ the composition condition is clearly sufficient and the
problem is to decide whether this condition is necessary.
Note also that since by L\"{u}roth theorem  (see e.g. \cite{sh}, p.13)
each field $k$ such that $\C \subset k \subset \C(z)$ and $k\neq \C$ 
is of the form $k=\C(R),$ $R\in \C(z)\setminus \C,$ 
it is easy to see that the conditions 
(*) hold if and only if the field $\C(P,Q)$ is a proper subfield
of $\C(z).$ In its turn, as $[\C(z):\C(P)]=\deg P,$ the last
condition is equivalent to the condition
$[\C(P,Q):\C(P)]<\deg P.$

Due to its connection with the center problem for the Abel equation
and with the classical Poincare center-focus problem for polynomial vector
field on the plane the polynomial moment problem has been studied 
in the recent papers \cite{c1}-\cite{ro}. In
particular, the truth of the composition conjecture was
established under additional assumption that $a,b$ are not critical
points of $P(z)$ (see \cite{c}) and under some other additional
assumptions (see \cite{c2},\cite{c3} and \cite{ro}). In this note we construct
a class of counterexamples to the composition conjecture.

\vskip 0.1cm

\noindent{\bf Claim 1.} {\it Let $B(z)$ and $D(z)$ be polynomials
such that $\C(B)\cap\C(D)$ contains a polynomial $P(z),$ $\deg P(z) >
1,$ and let $Q(z)=B(z)+D(z).$ Suppose that there exist $a,b \in \C$
satisfy $B(a)=B(b),$
$D(a)=D(b),$ $a\neq b.$ Then $m_i(P,Q,a,b)=0$ for all $i\geq 1.$} 

\vskip 0.1cm

\noindent{\it Proof.} Indeed, since $P(z)=A(B(z))$ for some polynomial
$A(z)$ and $B(a)=B(b)$ we have
$m_{i}(P,B,a,b)=m_{i}(A,z,B(a),B(b))=0$ for all $i\geq 1.$ Similarly,
$m_{i}(P,D,a,b)=0$ for all $i\geq 1.$ Therefore, also $m_{i}(P,Q,a,b)=0$
for all $i\geq 1.$ $\Box$
 
\vskip 0.1cm

\noindent{\bf Claim 2.} {\it Let $B(z)$ and $D(z)$ be polynomials
such that $\C(B)\cap\C(D)$ contains a polynomial $P(z),$ $\deg P(z) >
1,$ and let $Q(z)=B(z)+D(z).$ Suppose that $\C(B,D)=\C(z).$ Then 
$\C(P,Q)=\C(z)$.} 

\vskip 0.1cm

\noindent{\it Proof.} Assume the converse,
i.e. that the conditions (*) hold for
some polynomials $\tilde P(z),$ $\tilde Q(z),$ $W(z).$ 
To be definite suppose that $\deg B(z) \leq \deg D(z);$
then $\deg W(z) \vert \deg D(z).$
As $P(z)=C(D(z))$ for some polynomial $C(z),$ $\C(D)\cap
\C(W)$ contains a non-constant polynomial $P(z).$ By Engstrom theorem 
(see \cite{sh}, Theorem 5, p. 18) this condition follows that
$[\C(W,D):\C(D)]=\deg D/(\deg W,\deg D).$
Therefore, $(\deg W,\deg D)=\deg W > 1$ implies that
$[\C(D,W):\C(D)]<\deg D.$ Hence, $D(z)=\tilde D(F(z)), W(z)=\tilde
W(F(z))$ for some polynomials $\tilde D(z),$ $\tilde W(z),$
$F(z),$ $\deg F > 1.$ As $D(z)=\tilde D(F(z)),$ $Q(z)=\tilde
Q(W(z))=\tilde Q
(\tilde W(F(z)))$ and $B(z)=Q(z)-D(z)$ we see that $\C(B,D)\subset
\C(F).$ Since it contradicts to the condition $\C(B,D)=\C(z)$ we conclude
that $\C(P,Q)=\C(z)$. $\Box$

In order to get a counterexample $P(z), Q(z)$ to the composition
conjecture it is enough to find 
polynomials $B(z), D(z)$ which satisfy assumptions of both claims above.
Let us describe such polynomials.
Since $P\in \C(B)\cap \C(D),$ $\deg P >1,$ and $\C(B,D)=\C(z)$
the Engstrom theorem implies that $(\deg B,\deg D)=1.$ By 
the second Ritt theorem (see \cite{sh}, p. 24) 
conditions $P\in \C(B)\cap \C(D),$ $\deg P >1,$ and $(\deg B,\deg D)=1$
yield that up to linear change of variable either $B(z)= z^m,$ 
$D(z)= z^nR(z^m),$ where $R(z)\in \C[z]$ and $(n,m)=1$ (then 
$P(z)\in\C[z^{nm}R^m(z^m)]$) or $B(z)= T_n(z),$ 
$D(z)= T_m(z)$ for Chebyshev polynomials $T_n(z),$ $T_m(z)$ and $(n,m)=1$
(then $P\in\C[T_{nm}(z)]$).
In the first case, due to $(n,m)=1$, the conditions  
$B(a)=B(b),$
$D(a)=D(b),$ $a\neq b$ are equivalent to $a^m=b^m=\zeta,$ $a\neq b,$ 
where $\zeta$ is a root of $R(z).$ Therefore, whenever
$R(z)\neq z^l,$ $l\geq 1,$ we get a counterexample setting 
$Q(z)= z^m+z^nR(z^m),$ $P(z)\in\C[z^{nm}R^m(z^m)].$
In the second case, set 
$a=\alpha+1/\alpha,$ \linebreak $b=\beta+1/\beta,$ where $\alpha,$
$\beta$ satisfy $\alpha^m=\beta^m,$
$(\alpha\beta)^{n}=1,$ $\alpha\neq \beta, 1/\beta.$ 
Then $a\neq b$ and it follows from
$T_k(z+\frac{1}{z})=z^k+\frac{1}{z^k},$ $k \geq 1,$ that $T_n(a)=T_n(b),$
$T_m(a)=T_m(b).$ In this case counterexamples have a form:
$Q(z)=T_n(z)+T_m(z),$ $P(z)\in\C[T_{nm}(z)].$

\vskip 0.1cm

\noindent{\bf Acknowledgments.} I am grateful to P. M\"{u}ller, N.
Roytvarf and Y. Yomdin for interesting discussions. 
%Also I would like to
%thank the Weizmann Institute of Science and the Max Planck Institut
%f\"{u}r
%Mathematik for their support and hospitality.

\bibliographystyle{amsplain}

\end{document}